\documentclass[a4paper,10pt,twocolumn]{article}
\usepackage[utf8]{inputenc}

\usepackage[switch]{lineno}

\usepackage{hyperref}
\usepackage{bm}
\usepackage{authblk}
\usepackage{siunitx}
\usepackage{graphicx}
\usepackage{amsmath, amssymb}
\usepackage{subcaption}  
\usepackage{pgfplotstable} 

\usepackage[a4paper,left=1.5cm,right=1.5cm,top=2cm,bottom=2cm]{geometry}

\usepackage{nomencl}
\usepackage{ifthen}
\renewcommand{\nomgroup}[1]{%
\ifthenelse{\equal{#1}{C}}{\item[\textbf{Acronyms}]}{%
\ifthenelse{\equal{#1}{G}}{\item[\textbf{Greek Letters}]}{%
\ifthenelse{\equal{#1}{L}}{\item[\textbf{Latin Symbols}]}{%
\ifthenelse{\equal{#1}{S}}{\item[\textbf{Subscripts}]}{}}}}
}
\makenomenclature

\usepackage[numbers,sort&compress]{natbib}

\def\tsc#1{\csdef{#1}{\textsc{\lowercase{#1}}\xspace}}
\tsc{WGM}
\tsc{QE}
\tsc{EP}
\tsc{PMS}
\tsc{BEC}
\tsc{DE}

\begin{document}
\let\WriteBookmarks\relax
\def\floatpagepagefraction{1}
\def\textpagefraction{.001}

\title{A Predictive Surrogate Model for Heat Transfer of an Impinging Jet on a Concave Surface}
\author[1]{Sajad~Salavatidezfouli}
\author[2]{Saeid~Rakhsha}
\author[1]{Armin~Sheidani}
\author[3]{Giovanni~Stabile}
\author[1]{Gianluigi~Rozza}

\affil[1]{Mathematics Area, MathLab, International School for Advanced Studies (SISSA), Trieste, Italy}

\affil[2]{Department of Mechanical Engineering, Semnan University, Semnan, Iran}

\affil[3]{Department of Pure and Applied Sciences, Informatics and Mathematics Section, University of Urbino Carlo Bo, Urbino, Italy}

\date{} 

\twocolumn[
  \begin{@twocolumnfalse}
    \maketitle
	   \begin{abstract}
            This paper aims to comprehensively investigate the efficacy of various Model Order Reduction (MOR) and deep learning techniques in predicting heat transfer in a pulsed jet impinging on a concave surface. Expanding on the previous experimental and numerical research involving pulsed circular jets \citep{rakhsha2023effect}, this investigation extends to evaluate Predictive Surrogate Models (PSM) for heat transfer across various jet characteristics. To this end, this work introduces two predictive approaches, one employing a Fast Fourier Transformation augmented Artificial Neural Network (FFT-ANN) for predicting the average Nusselt number under constant-frequency scenarios. Moreover, the investigation introduces the Proper Orthogonal Decomposition and Long Short-Term Memory (POD-LSTM) approach for random-frequency impingement jets. The POD-LSTM method proves to be a robust solution for predicting the local heat transfer rate under random-frequency impingement scenarios, capturing both the trend and value of temporal modes. The comparison of these approaches highlights the versatility and efficacy of advanced machine learning techniques in modelling complex heat transfer phenomena.


  \end{abstract}
  \end{@twocolumnfalse}

]
\maketitle
\printnomenclature

\section{Introduction}
Owing to their high thermal efficiencies, impinging jets encompass a wide range of industrial applications including cooling of turbine blades, electronic devices, combustors walls, heat exchangers and so on \citep{raizner2020effect, selimefendigil2014pulsating, sarkar2023review, kim2023computational, kim2022enhancing}. Considering the applications which are inherently associated with a high level of oscillation like the gas turbine blades, cardinal importance is placed on the numerical and experimental study of pulsating impinging jets which has attracted a great deal of attention in recent years \citep{ hadipour2018heat, van2018experimental, attalla2020experimental, ikhlaq2021flow, parida2021heat}. It should be noted that one of the main objectives of studying the complicated physics of such jets is to devise a way to improve cooling efficiency. In this regard, the most common strategies to reach this goal have been inlet perturbation, jet pulsation, utilization of the rough surface, and change of the nozzle geometry. \citep{rakhsha2023effect, selimefendigil2022combined, dirker2021influence, zhao2023numerical, kumar2023experimental, cheng2023pore, li2023random}. 

The improvement strategies can be generally divided into active/passive methods. For instance, alterations to parameters such as the nozzle diameter or the distance from the surface could be regarded as examples of passive methods, while techniques dealing with the velocity or the frequency of the impinging jet are considered as an active strategy. Therefore, in order to improve the heat transfer for an impingement jet, all the aforementioned parameters might be taken into consideration \citep{hussain1989elliptic, husain1991elliptic, husain1993elliptic, mi2000centreline, singh2003mixing}. 

To illustrate, \cite{alimohammadi2015numerical} in a numerical and experimental study proposed a correlation describing the relation between the Nusselt and the Reynolds number with the distance between the nozzle and the surface. As for the importance of the nozzle geometry, \cite{lee2000effect} showed that the utilization of elliptical jets results in a higher thermal efficiency compared to their circular counterparts. In this regard, \cite{koseoglu2010role} reported that the increment in the nozzle aspect ratio is directly proportional to the growth in the average Nu number on a flat plate. Owing to the importance of this issue, the investigation of the effect of the nozzle design parameters on the heat transfer of impinging jets has continued in recent years. To illustrate, \cite{paolillo2022effects} investigated the influence of swirl and Re number along with the distance of a nozzle on the heat transfer. The results showed that the nozzle distance significantly affects the heat transfer and as a result, it was suggested to utilize the jet at a close distance along with swirl addition for higher thermal efficiency. \cite{khan2023comparison}, in a numerical study, drew a comparison between the oscillating (sweeping) and constant jets for cooling the blade of a gas turbine. The results proved that the employment of sweeping jets leads to a considerable improvement in heat transfer. \cite{mohammadshahi2021experimental} in an experimental study, compared the flow characteristics and heat transfer of an oscillating and steady jet in a cross flow. It was reported that the heat transfer is mainly affected by the cross-flow in the case of the oscillating jet causing the alteration of the Blowing Ratio (BR), while as for the constant jets the rise in BR results in a significant increase in the heat transfer. \cite{feng2021experimental} conducted a numerical and experimental study on the effect of various jet parameters, including the jet shape, on heat transfer. The upshot of their findings was that circular jets outperformed fan jets due to their higher heat transfer coefficient.

Due to the complexities associated with the physics of heat transfer in industrial applications, Machine Learning (ML) techniques have been employed to facilitate the understanding of such phenomena in recent years \citep{xie2022two, lee2022classification, hachem2021deep, bhattacharyya2021application}. To illustrate, \cite{mohammadpour2022machine} implemented various ML techniques including k-nearest neighbour (k-NN), Random Forest (RF) and Multi-Layer Perceptron (MLP) to assess their performance regarding the numerical simulation of nano-fluid heat transfer in a micro-channel equipped with double synthetic jets. The results revealed that the k-NN algorithm yields more accurate predictions compared to the other methods employed. It should be noted that in the work of \cite{zhou2020machine} a comparison between the ML techniques including Artificial Neural Networks (ANN), RF, AdaBoost and Extreme Gradient Boosting (XGBoost) suggested that the ANN outperforms the other methods for the prediction of condensation coefficient in micro-channels. \cite{zhu2021machine} trained an ANN based on experimental data to predict the boiling and condensation for a refrigerant in a channel proving that the use of parameterized ANN yields accurate results. \cite{oka2020parameter} employed Kalman Filter (KF) to estimate the heat flux in casting molds. The results suggested that KF offers a great performance for the estimation of heat flux in Inverse Heat Transfer Problems (IHTP). As for the applications of ML in IHTPs, \cite{zhu2022deep} proposed a high-precision hybrid deep learning model based on Conventional Neural Networks (CNN) and Long Short-Term Memory (LSTM) for the time-dependent parameters. \cite{kianimoqadam2023calculating} predicted the radiation view factors for particles and face neighbours utilising LSTM and Gated Recurrent Unit (GRU).
\\
Based on the studies conducted thus far, it is evident that research on impinging jets has been primarily focused on exploring various design parameters and their impact on thermal efficiency. In other words, the implementation of modern techniques, specifically ML, to mitigate the computational costs associated with the numerical simulation of impinging jets has received limited attention in this field. To address this issue, the primary objective of this study is to utilize a range of deep learning techniques, such as frequency-based MLP, LSTM, and Proper Orthogonal Decomposition (POD)-Galerkin, for predicting heat transfer. Specifically, this research aims to develop a novel deep-learning model for predicting the Nu number on a concave plate subjected to a pulsating impinging jet. Fundamentally, the prediction of Full Order Modelling (FOM) by ML methods is associated with difficulties mainly due to the presence of noise. In this regard, this study is intended to bring this issue into sharp focus. The model will consider a wide range of design parameters, including inlet velocity, frequency, and nozzle-to-surface distance as to the best of the authors' knowledge there is yet to be any report on the prediction of the flow and geometrical parameters. In the upcoming section, this paper begins by providing a detailed problem description. Subsequently, the numerical method, encompassing the governing equation, turbulence model, geometry and grid, and validation, is presented. The subsequent sections elaborate on the deep learning models and nonlinear reduced-order models. Finally, the paper concludes with the presentation of results and discussions, followed by the conclusion.          

\section{Full-Order Model}
The numerical simulation of the heat transfer on the concave surface has been carried out with ANSYS FLUENT V2019 \citep{fluent2019ansys}. In the following section, the details of the full-order model are presented. 

\subsection{Problem Description}
\begin{figure}
	\centering
	\includegraphics[scale=.6]{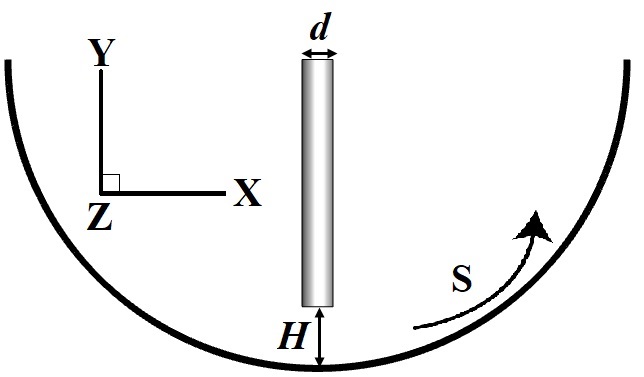}
	\caption{Schematic of the computational model}
	\label{FIG_schematic}
\end{figure}

Fig. \ref{FIG_schematic} demonstrates the schematic of the computational domain consisting of a concave surface under constant heat flux prone to a cooling jet. This model is selected based on the study of \cite{rakhsha2023effect}, where the effects of various nozzle shapes (circle, rectangle, elliptic, and square) and nozzle-to-surface distance, $ H/d$ were investigated.
Given our main objective, which is the prediction of heat transfer rate on the surface, we concentrate solely on the optimal nozzle shape, i.e., the circle \citep{rakhsha2023effect} accompanied with variable nozzle-to-surface distance. The concave surface is semi-cylindrical with a diameter of $240 \si{mm}$ and a length of $120 \si{mm}$. The nozzle diameter of the circular jet was $12 \si{mm}$. 

Additionally, we explore two different approaches for jet velocity; harmonic velocity and random inlet velocity. These approaches will be discussed in detail in the next sections.

\subsection{Governing Equation}
The 3-dimensional numerical simulation is carried out to investigate the effect of the geometrical and flow parameters on the heat transfer of the concave surface. In this section, the general properties of the fluid, governing equations, and solver properties will be discussed. Moreover, a further discussion on turbulence modelling will be presented.

As for the fluid properties, air with constant thermophysical properties is considered. The fluid properties are presented in Table \ref{table1}. 

\begin{table*}[]
\caption{Thermophysical properties of the fluid}
\label{table1}
\begin{tabular}{cccccccc}
$\rho$                     & $\mu$       & $k$     & $C_p$   \\[1ex] \hline
($kg/m^3$) & ($Pa.s$)     & ($W/mK$)  & ($J/kg.K$) \\[1ex] \hline
1.225                   & 1.789e-5 & 0.024 & 1006  
\end{tabular}
\end{table*}

The governing equations of the problem including Navier-stokes and energy equations are as follows:
\begin{equation}
\begin{aligned}
&-\frac{\partial \bar{u}_i}{\partial x_i}=0,\\
\end{aligned}
\end{equation}
\begin{equation}
\begin{aligned}
&\rho \frac{\partial u_i}{\partial t}+\frac{\partial\left(\overline{u_i} \overline{u_j}\right)}{\partial x_j}=-\frac{\partial \bar{P}}{\partial x_j}+\frac{\partial}{\partial x_j}\left[\mu\left(\frac{\partial \overline{u_i}}{\partial x_j}+\frac{\partial \overline{u_j}}{\partial x_i}\right)-\rho \overline{u_i u_j}\right],\\
\end{aligned}
\label{eq4}
\end{equation}
\begin{equation}
\begin{aligned}
&\frac{\partial \bar{T}}{\partial t}+\frac{\partial\left(u_i \bar{T}\right)}{\partial x_i}=\frac{\partial}{\partial x_i}\left(\alpha \frac{\partial \bar{T}}{\partial x_i}-\overline{u_i T}\right).\\
\label{eq5}
\end{aligned}
\end{equation}

\subsection{Turbulence Model}
Regarding equations \ref{eq4} and \ref{eq5}, $\rho \overline{u_i u_j}$ and $\overline{u_i T}$ are the Reynolds stress tensor and turbulent heat flux, respectively. These terms can be modelled by the well-known 1- or 2-equation turbulence models as follows:

\begin{equation}
\begin{aligned}
&\overline{\rho u_i^{\prime} u_j^{\prime}}=\mu_\psi\left(\frac{\partial u_i^{\prime}}{\partial x_j}+\frac{\partial u_j^{\prime}}{\partial x_i}\right)-\frac{2}{3} \delta_{i j} \rho k,\\
\end{aligned}
\end{equation}
\begin{equation}
\begin{aligned}
&\rho \overline{u_{\imath}^{\prime} T^{\prime}}=\frac{\mu_t}{\operatorname{Pr}}\left(\frac{\partial T}{\partial x_j}\right),
\end{aligned}
\end{equation}

where $\mu_t$ is turbulent viscosity, $\delta_{i j}$ is Kronecker delta and $k = u_i u_j$ is turbulent kinetic energy. To determine $\mu_t$ the $k-\omega$ turbulence model was employed which can rather precisely address the near wall phenomena. The $k-\varepsilon$ has been reported previously to outperform other turbulence models in terms of impingement jets, which will be discussed in section \ref{validation study}.

\subsection{Computational Grid}
\begin{figure}
	\centering
	\includegraphics[scale=0.5]{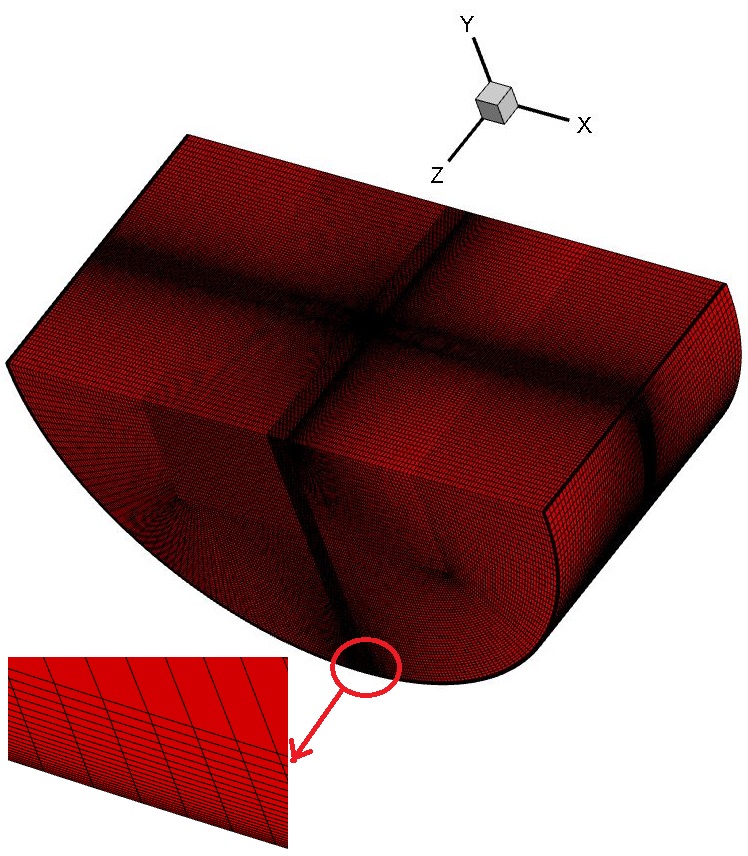}
	\caption{Structured grid for the model with a closer look at the boundary layer mesh near the concave surface}
	\label{FIG_Grid}
\end{figure}

\begin{figure}
	\centering
	\includegraphics[scale=0.475]{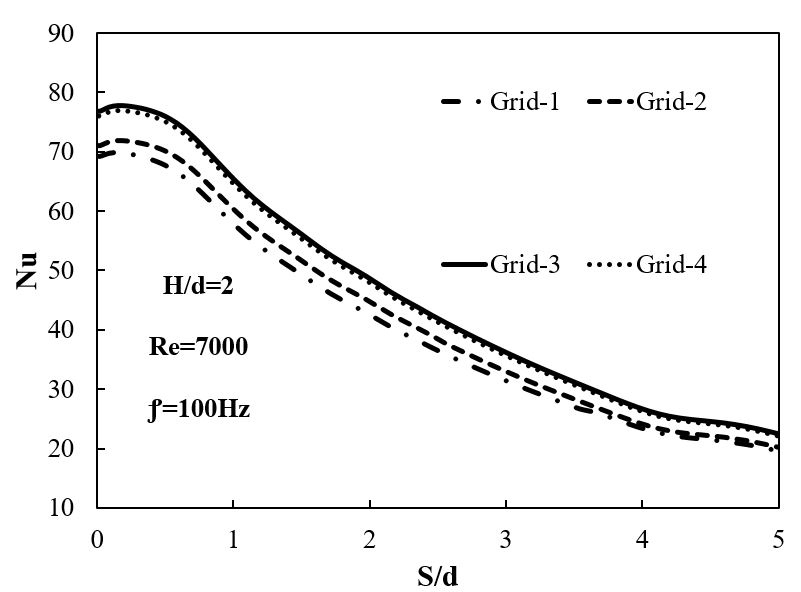}
	\caption{Variation of local Nusselt number for different meshes}
	\label{FIG_GridStudy}
\end{figure}

A fully structured mesh has been considered for the simulation as 
shown in Fig. \ref{FIG_Grid}. To examine the independence of the numerical solution to the element size, four sets of meshes with respectively $4\times10^5$, $8\times10^5$, $1.5\times10^6$, and $2.5\times10^6$ is generated.
The local Nusselt number demonstrates represents the dimensionless heat transfer coefficient on a surface experiencing convective heat transfer. Fig. \ref{FIG_GridStudy} demonstrates the variation of local Nusselt number along the concave surface (S-direction based on Fig. \ref{FIG_schematic}) for different grids. Clearly, the third grid corresponding to 1500000 hexahedral elements has been selected for other simulations. For all the geometries (in terms of different H/d) a fine mesh was generated near the concave surface to achieve the best accuracy of the numerical solution. It is worth noting that for capturing the near-wall physics, an element height of $10^{-6} \si{m}$ corresponding to $Y^+<1$ is employed for the concave surface. 

\subsection{Boundary and Initial Condition}
\begin{figure}
	\centering
	\includegraphics[scale=0.475]{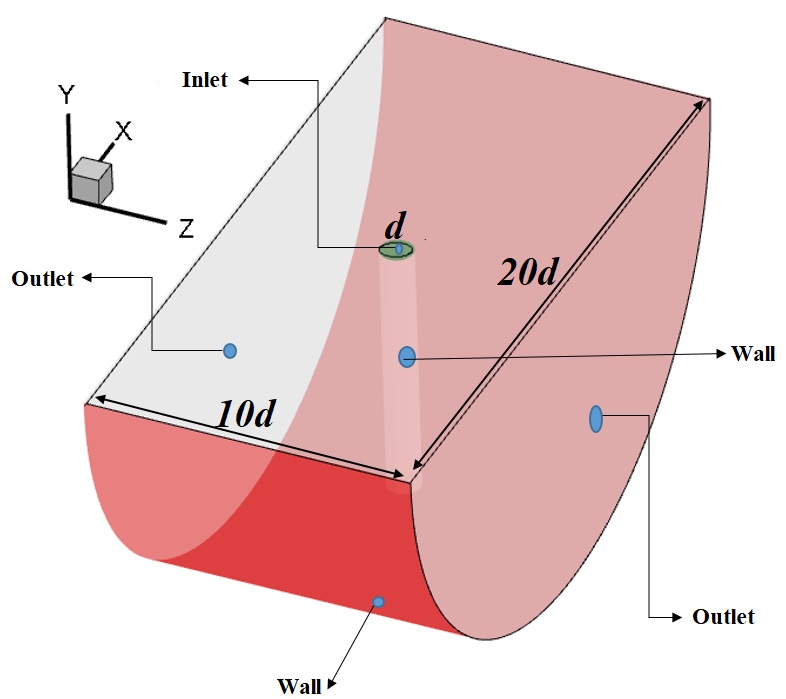}
	\caption{Utilized boundary conditions for the CFD simulation}
	\label{FIG_BC}
\end{figure}

The boundary conditions of the CFD model are shown in Fig. \ref{FIG_BC}. The inlet temperature was set to $300 \si{K}$ with harmonic/random velocity based on the prediction purpose. No-slip assumption along with constant heat flux was applied for the concave surface. Details of boundary conditions have been presented in Table \ref{table_bc}. Air was modelled as incompressible flow since the Mach number is far less than 0.3 

\begin{table}[]
\caption{Values of the boundary conditions}
\label{table_bc}
\begin{tabular}{cc}
\hline
Inlet           & \begin{tabular}[c]{@{}c@{}}Velocity inlet: Harmonic/Random \\ Inlet temperature: $300 \si{K}$ \\ Turbulence intensity: 5\% \\ Turbulence kinetic energy : 0.24\si{m^2/s^2} \end{tabular}    \\[1ex] \hline
Outlet          & Pressure Outlet:   Gauge pressure : 0 Pa     \\[1ex] \hline
Concave & \begin{tabular}[c]{@{}c@{}}Constant heat flux-$2500 \si{W/m^2}$\\ No-slip condition \end{tabular}                                               \\ \hline
\end{tabular}
\end{table}

\subsection{Solver Properties}
The finite volume method (FVM) was applied to discretize the governing equations. For all the equations the second-order upwind scheme was implemented to address advection terms. SIMPLE-C algorithm has been used for pressure-velocity coupling. Moreover, the time discretization has been done by the second-order implicit scheme. The timestep of the transient simulation was determined by dividing the jet pulse period by 100. To ensure the periodicity of the final solution, the simulation continued until 20 cycles. Convergence criteria were set to $10^{-4}$ for all equations in each time step.

\subsection{Validation Study}\label{validation study}
\begin{figure}
	\centering
	\includegraphics[scale=0.44]{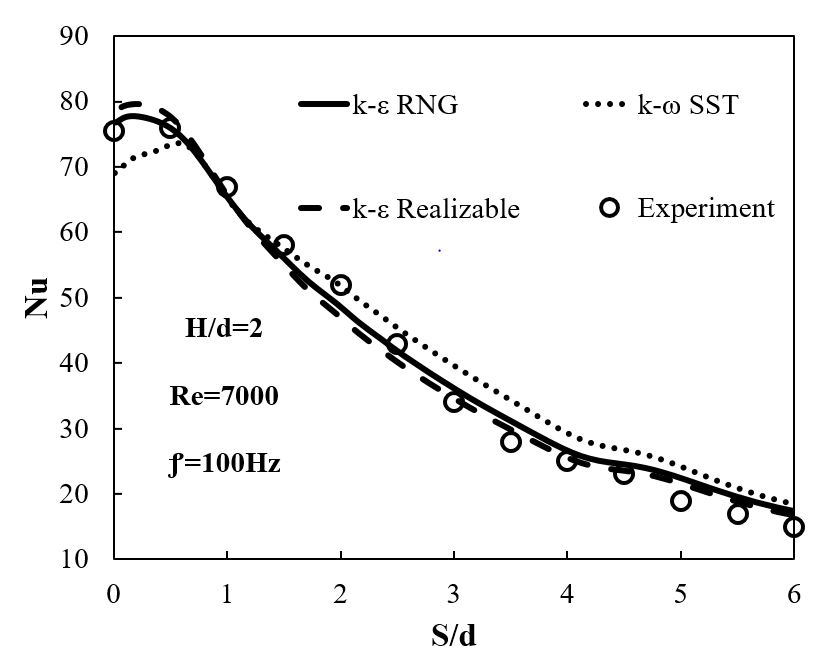}
	\caption{The performance of turbulence models to predict the local Nu number distribution}
	\label{FIG_Validation}
\end{figure}

The local Nusselt number based on the available experimental data and numerical simulation of the concave surface along the S-direction is shown in Fig. \ref{FIG_Validation}. A good agreement was observed between all turbulence models and the experimental data of \cite{rakhsha2023effect}. However, they observed a comparable error for all turbulence models in the stagnation region, i.e. $S/d<1$ and $X/d<1$. This is due to the fact that as the jet impinges the surface, a stagnation region is generated corresponding to the values of low Reynolds numbers. 

\subsection{Physics of the Flow}
\begin{figure}
	\centering
	\includegraphics[scale=.22]{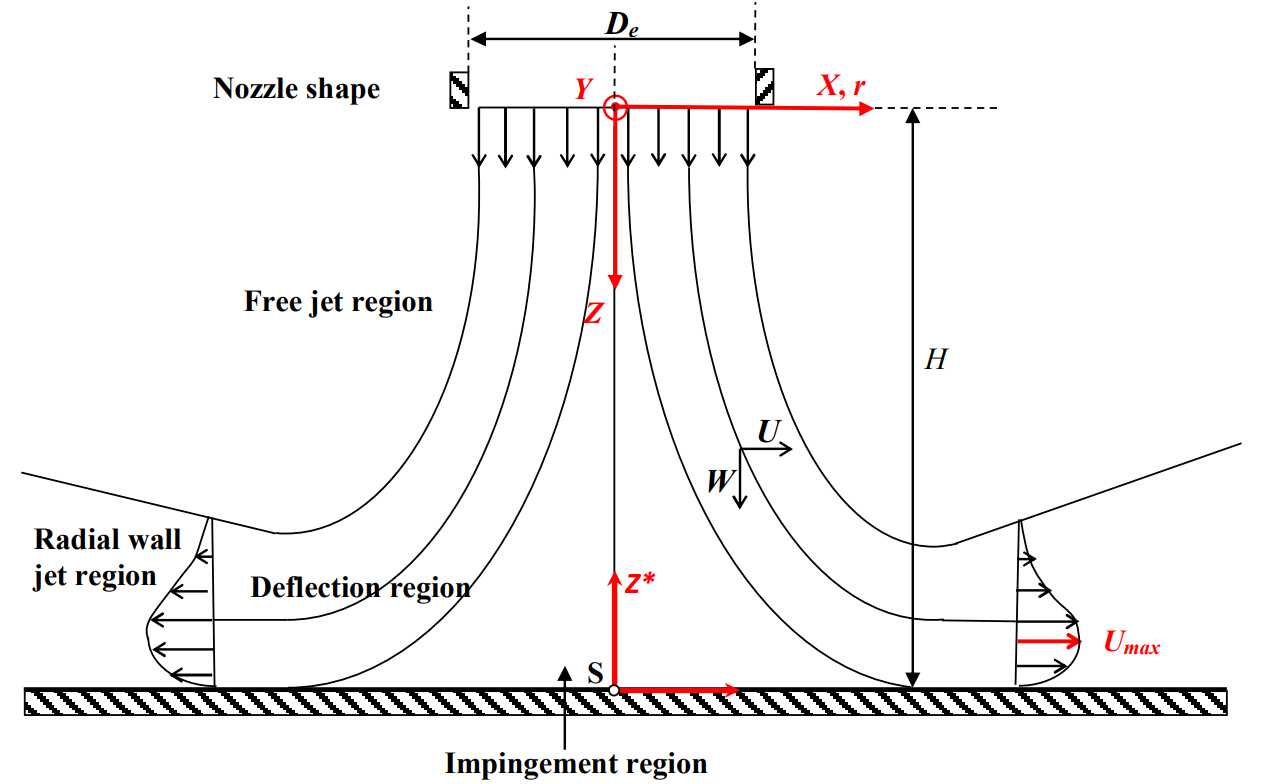}
	\caption{Formation of different regions in the context of the impinging jet on a plate \cite{sodjavi2015impinging} }
	\label{FIG_flowPhysics}
\end{figure}

Before exploring predictive models, we need to understand the physics of the flow in terms of impingement cooling. According to Fig. \ref{FIG_flowPhysics}, the flow field of an impinging jet consists of three regions. The free jet, the impingement region and the wall jet region. These regions are highly sensitive to the effective flow and geometrical parameters including Reynolds number, nozzle-to-plate distance, nozzle section shape and surface geometry. Change in any of the mentioned parameters results in changes in the boundaries between three regions, which also subsequently affects the heat transfer from the impinging surface. Among all the available parameters, the jet velocity can be considered as the most potential variable in active flow control of the impinging jet. A change in the jet velocity and frequency causes a change in the vortex shedding frequency in the flow field, which eventually changes the formation and thickness of the boundary layer in the region of the wall jet.

Moreover, in practical cases, it is difficult to find a general correlation that can provide an explicit relationship between the heat transfer coefficient or surface temperature in terms of jet velocity. Hence, the use of advanced ML techniques such as DRL is necessary for such systems.

\section{Deep Learning Models} \label{dlm}
In the following section, we explore three key architectures—Multilayer Perceptrons (MLP), Long Short-Term Memory (LSTM) and Transformers networks. These models offer unique capabilities, from MLP's versatility in capturing non-linear relationships to LSTM's proficiency in handling sequential data, and Transformer's focus on long-range dependencies.

\subsection{Feed-Forward Network}
\begin{figure*}
	\centering
	\includegraphics[scale=0.3]{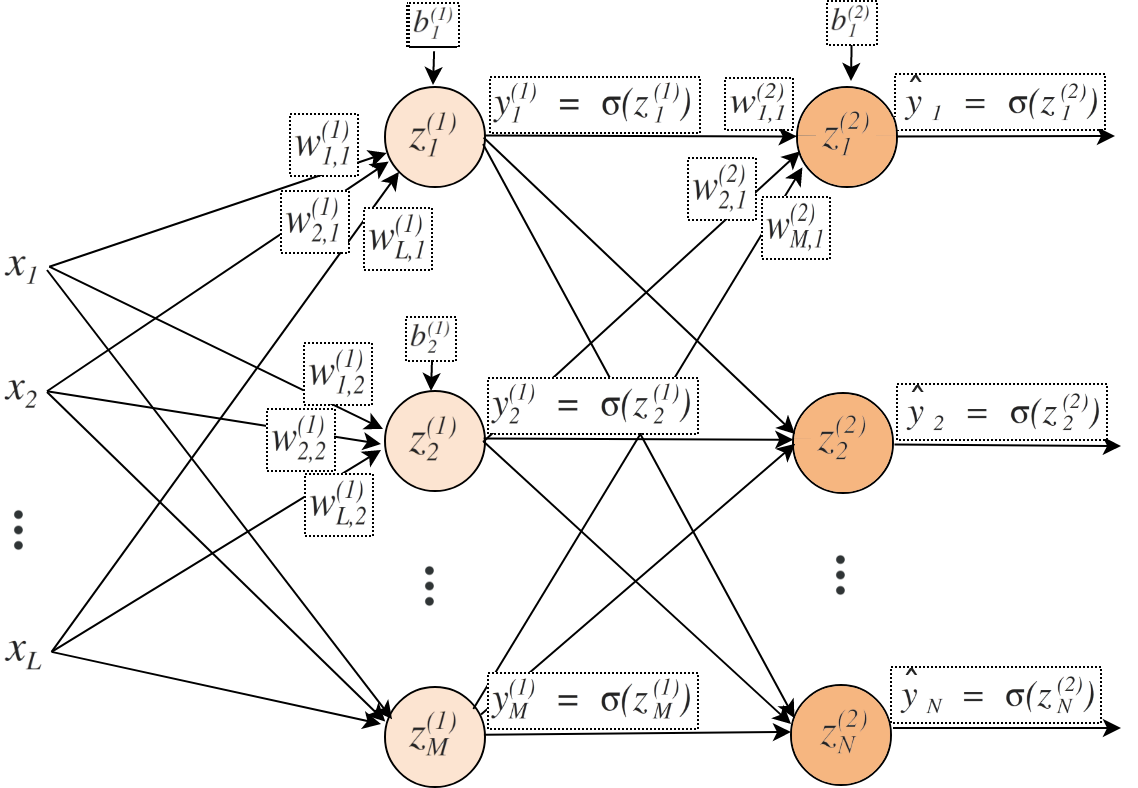}
	\caption{Structure of the feed-forward network}
	\label{FIG_MLPNet}
\end{figure*}

To generate functional relationships effectively using the deep-learning platform between the input and output, the first simplest method is the feed-forward network, also known as multilayer perception. The schematic view of MLP is shown in Fig. \ref{FIG_MLPNet}. As shown in the figure, the input signals are denoted by the array $x = [x_1, x_2, x_3, \ldots, x_N]$. As these signals are input to the neuron (shown in coloured circles), they undergo multiplication by their corresponding synaptic weights, forming the array $w = [w_1, w_2, w_3, \ldots, w_N]$. This multiplication produces the value $z$, referred to as \textit{activation potential}, as defined by the following formula:

\begin{equation}
\label{eq_MLP}
\begin{aligned}
z = \sum_{i=1}^{N} x_i \cdot w_i + b,
\end{aligned}
\end{equation}
where, the additional term $b$, unaffected by the input array, is the bias, introducing an extra degree of freedom to the model. Subsequently, the value $z$ undergoes transformation through an activation function $ReLU$. This function serves to constrain the value within a specific range, ultimately yielding the neuron's final output $y$.

\subsection{Long Short Term Memory}
\begin{figure}
	\centering
	\includegraphics[scale=0.0275]{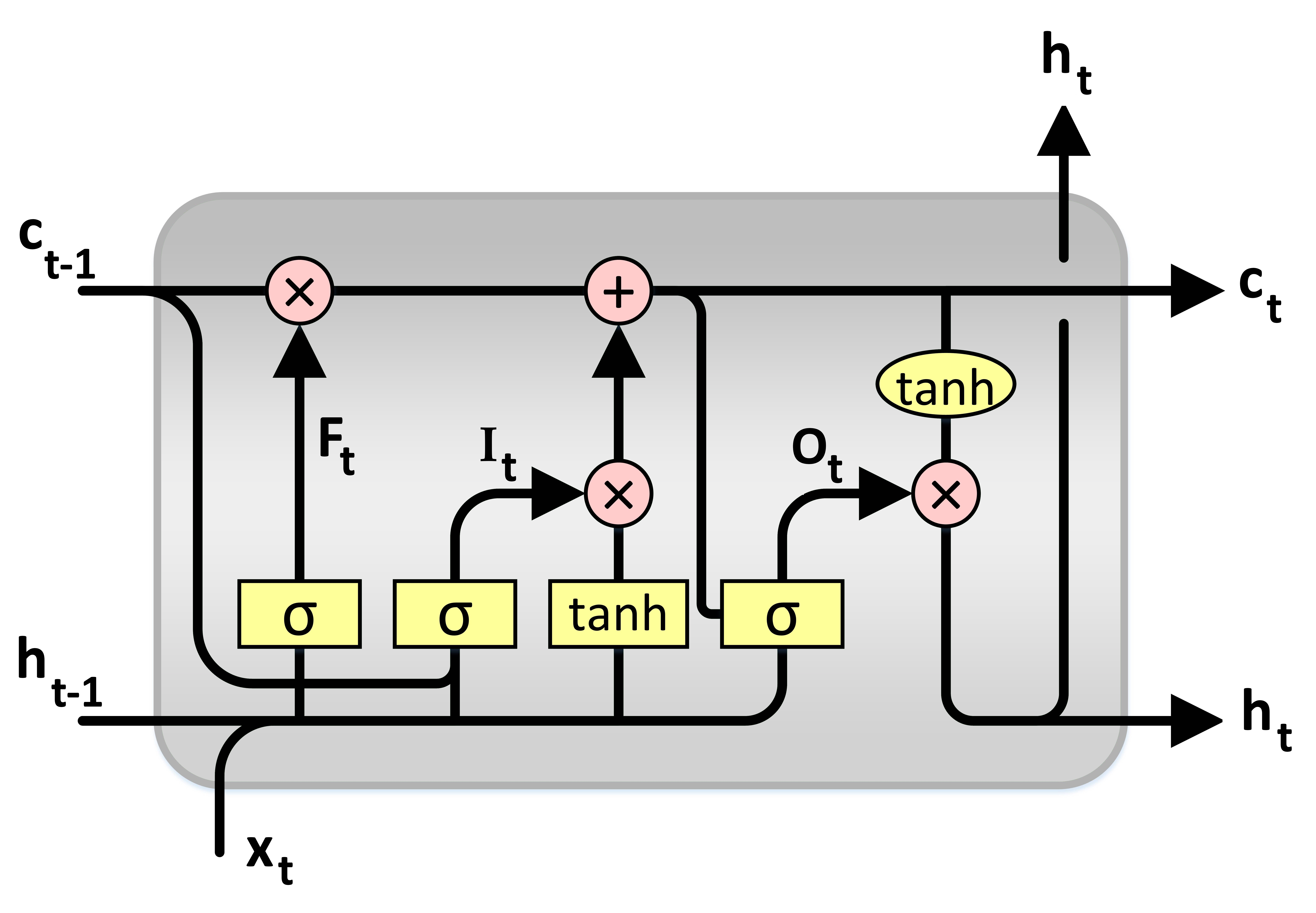}
	\caption{Structure of the LSTM block}
	\label{FIG_LSTMnet}
\end{figure}

The idea of LSTM networks was first proposed by \cite{hochreiter1997long} which was developed in the later studies \cite{van2020review}. In Fig. \ref{FIG_LSTMnet} the schematic architecture of LSTM has been presented. 

One of the main advantages of LSTMs over their classic counterparts known as Recurrent Neural Networks (RNN) is their ability to recall a longer span in the time-series. More precisely, as seen in Fig. \ref{FIG_LSTMnet}, an LSTM network enjoys two inputs. The input denoted by C stands for "Cell State" which is another term for the long-term memory that RNNs lack. The cell state connects all the stages of the network together going through the forget and input gates. The formulation of these gates and cell states is expressed as follows:

\begin{equation}
\label{eq_LSTM1}
\begin{aligned}
& I_t=\sigma\left(W_{x I} x_t+W_{h I} h_{t-1}+W_{c I} c_{t-1}+B_I\right),
\end{aligned}
\end{equation}
\begin{equation}
\label{eq_LSTM2}
\begin{aligned}
& F_t=\sigma\left(W_{x F} x_t+W_{h F} h_{t-1}+W_{c F} c_{t-1}+B_F\right), 
\end{aligned}
\end{equation}
\begin{equation}
\label{eq_LSTM3}
\begin{aligned}
& c_t=F_t c_{t-1}+I_t \tanh \left(W_{x c} x_t+W_{h c} h_{t-1}+B_c\right),
\end{aligned}
\end{equation}
\begin{equation}
\label{eq_LSTM4}
\begin{aligned}
& O_t=\sigma\left(W_{x O} x_t+W_{h o} h_{t-1}+W_{c o} c_t+B_O\right),
\end{aligned}
\end{equation}
\begin{equation}
\label{eq_LSTM5}
\begin{aligned}
& h_t=o_t \tanh \left(c_t\right),
\end{aligned}
\end{equation}
where, $W_{xI}$, $W_{xF}$, $W_{xc}$, $W_{xO}$ represents matrices of weight from the input gate, forget gate, long-term cell state, output gate to input $x_t$, respectively. $W_{hI}$, $W_{hF}$, $W_{hc}$, $W_{hO}$ represents matrices of weight from the input gate, forget gate, long-term cell state, output gate to intermediate output $h_{t-1}$, respectively. $W_{cI}$, $W_{cF}$, $W_{cO}$ represents matrices of weight from input gate, forget gate, output gate to cell state $c_t$, respectively \citep{dasan2021novel}. $B$ is the Bias vectors and $\sigma (.)$ is the sigmoid function.

\subsection{Transformer}
\begin{figure}
	\centering
	\includegraphics[scale=0.11]{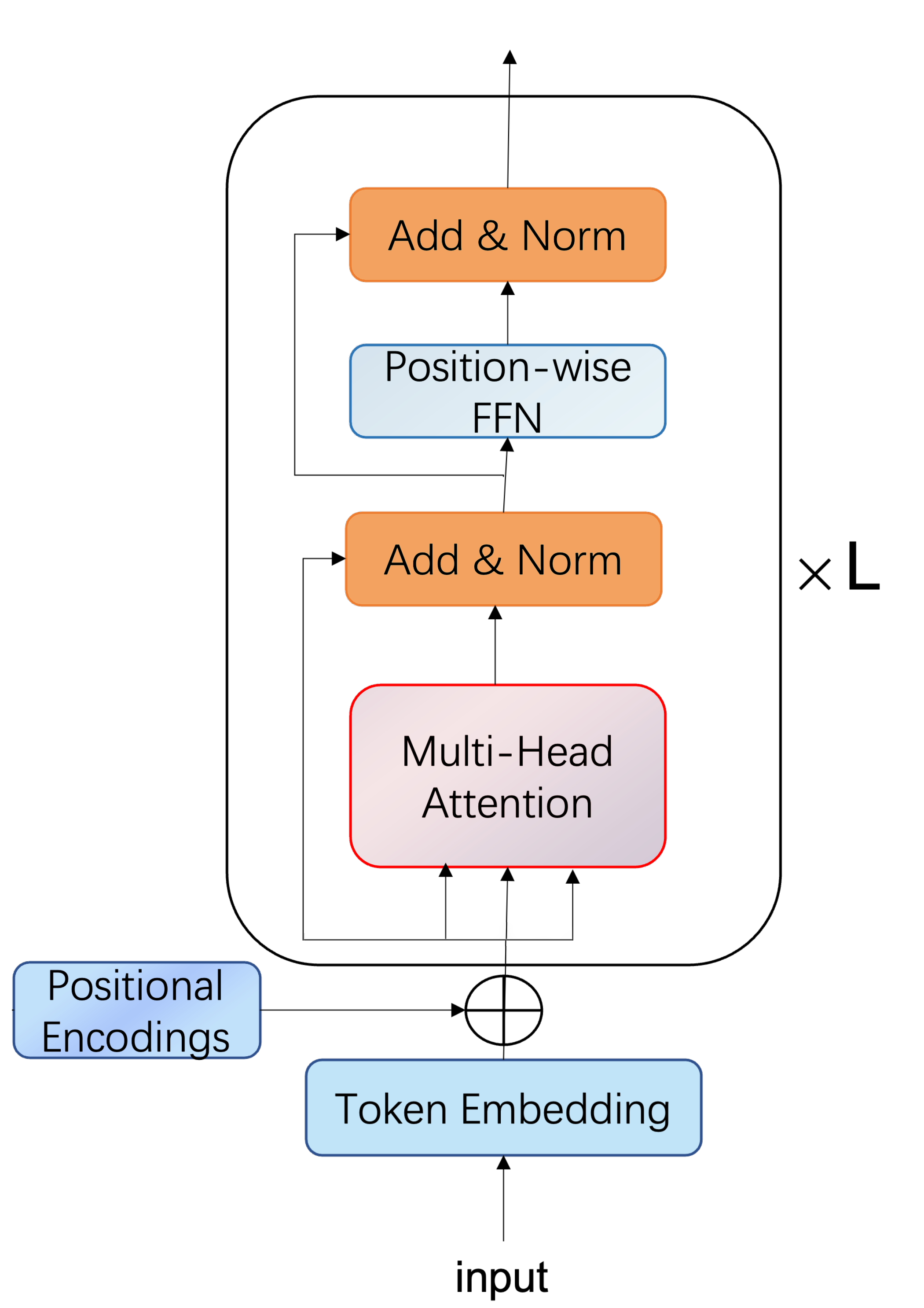}
	\caption{Structure of the encoder-only transformers network}
	\label{FIG_transformersnet}
\end{figure}

The concept of Transformer models goes back to the study "Attention is all you need" introduced by \cite{vaswani2017attention} which revolutionized the field of Natural Language Processing (NLP). Owing to their great performance, Transformers have been utilized in many applications including image processing techniques, denoising and machine translation. The architecture of Transformers relies heavily on the concept of "self-attention mechanisms" at the core. The self-attention mechanism enables the model to establish a relation between the relevant data points in the time series. In Fig. \ref{FIG_transformersnet} a schematic overview of the \textit{encoder-only} Transformer (also known as \textit{vanilla} transformer) network has been presented. 

Unlike LSTM, the Transformer has no recurrence. Instead, it utilizes the positional encoding added in the input embeddings, to maintain the sequence information. The following explanation explains the details of the transformer network.

\textbf{Positional Encoding}: In an encoder-only Transformer, the encoding vector for each position index \(t\) is defined as follows:
\begin{equation}
\label{eq_trans1}
\begin{aligned}
   PE(t)_i = 
   \begin{cases} 
   \sin(\omega_i t), & \text{if } i \mod 2 = 0, \\
   \cos(\omega_i t), & \text{if } i \mod 2 = 1,
   \end{cases}
\end{aligned}
\end{equation}
where $\omega_i$ represents the frequency for each dimension.

\textbf{Multi-head Attention}: The scaled dot-product attention used by the Transformer with the Query-Key-Value (QKV) model is expressed as:
\begin{equation}
\label{eq_trans2}
\begin{aligned}
\text{Attention}(Q, K, V) = \text{softmax}\left(\frac{K^T}{\sqrt{D_k}}\right) V,
\end{aligned}
\end{equation}
where queries $(Q \in \mathbb{R}^{N \times D_k})$, keys $(K \in \mathbb{R}^{M \times D_k})$, values $(V \in \mathbb{R}^{M \times D_v})$, and $(N, M)$ denote the lengths of queries and keys (or values). $D_k$ and $D_v$ represent the dimensions of keys (or queries) and values, respectively. The Transformer employs multi-head attention with $H$ different sets of learned projections instead of a single attention function.

\textbf{Feed-forward and Residual Network}:
The feed-forward network is similar to that explained eariler. In a deeper module, a residual connection module followed by a layer normalization module is inserted around each module:

\begin{equation}
\label{eq_trans3}
\begin{aligned}
H_0 = \text{LayerNorm}(\text{SelfAttn}(X) + X), \\
H = \text{LayerNorm}(FFN(H_0) + H_0),
\end{aligned}
\end{equation}
where $(\text{SelfAttn}(\cdot))$ denotes the self-attention module, and $(\text{LayerNorm}(\cdot))$ denotes the layer normalization operation.

\section{Reduced-Order Model}
The POD is the most widely used technique to reduce data and extract an optimal set of orthonormal bases. In this formulation, basis functions are constructed by collecting temporal snapshots for the variable of interest during the full model solution. Each snapshot vector, $\tilde{Nu}$, holds the values of the local Nu number at the center of the computational cell. The sampled values at the snapshot $s$ are stored in the matrix $\mathbf{\tilde{Nu}}_{s}$ with $\mathcal{N}$ entries ($\mathcal{N}$ being the number of nodes) to construct the snapshot matrices of $\mathbf{\tilde{Nu}}=\left({\tilde{Nu}}_1, \ldots, {\tilde{Nu}}_s, \ldots, {\tilde{Nu}}_S\right)$. For the sake of simplicity, the details will be provided for a general snapshot matrix $\varphi$.

The goal of POD is to find a set of orthogonal basis functions $\left\{\phi_s\right\}, s \in\{1,2, \ldots, S\}$, such that it maximizes:
\begin{equation}
\label{eq_pod1}
\begin{aligned}
& \frac{1}{S} \sum_{s=1}^S\left|<\varphi, \phi_s>_{L^2}\right|^2,
\end{aligned}
\end{equation}
subject to:

\begin{equation}
\label{eq_pod2}
\left\langle\phi_i, \phi_j\right\rangle_{L_2}=\delta_{i j}
\end{equation}
where $\langle\cdot, \cdot\rangle_{L^2} $ is the canonical inner product in L2 norm.

The approach introduced by \cite{sirovich1987turbulence} seeks to find an optimal set of basis functions $\phi$ for the optimization problem, i.e. Eq. \eqref{eq_pod1}. This requires performing a Singular Value Decomposition (SVD) of the snapshot matrix $\varphi$ given in the form,

\begin{equation}
\label{eq_pod3}
\begin{aligned}
& \varphi = U \Sigma V^T.
\end{aligned}
\end{equation}

The terms $U \in R^{\mathcal{N} \times \mathcal{N}}$ and $V \in R^{\mathcal{S} \times \mathcal{S}}$ are the matrices that consist of the orthogonal vectors for $\varphi \varphi^{T}$ and $\varphi^{T}\varphi$, respectively and $\Sigma$ is a diagonal matrix of size $\mathcal{N} \times \mathcal{S}$. The non-zero values of $\Sigma$ are the singular values of $\varphi$, and these are assumed to be listed in order of their decreasing magnitude \citep{xiao2015non}. Further on, the POD vectors can be defined as the column vectors of the matrix $U$. These vectors are considered to be optimal in the sense that no other rank $N$ set of basis vectors can be closer to the snapshot matrix $\varphi$ in the Frobenius norm. 

A reduced-order approximation of the field can be described as follows:

\begin{equation}
\label{eq_pod4_5}
\varphi \approx \sum_{n=1}^N \alpha_n(t) \cdot \phi_n.
\end{equation}

The loss of information due to the truncation of the POD expansion set to $N$ vectors can be quantified as follows:

\begin{equation}
\label{eq_pod5}
\begin{aligned}
& L=\frac{\sum_{n=1}^N \lambda_n^2}{\sum_{n=1}^S \lambda_n^2},
\end{aligned}
\end{equation}
where $\lambda$ denotes singular values.

\section{Results and Discussion}
In this section, the results obtained by employing different model order reduction/deep learning techniques, namely MLP, LSTM, Transformer, and the combination of POD and LSTM, for the two scenarios will be investigated. These scenarios consist of a constant-frequency and a multiple-frequency impingement jet. The objective is to analyze and compare the performance of these techniques in each of the scenarios and draw meaningful conclusions from the results.

\subsection{Constant-Frequency Impingement Jet}
As for the simplest scenario, it was assumed that the outlet velocity of the jet oscillated with time sinusoidally with one frequency. The frequency value ranged from $5 Hz$ to $100 Hz$ for different cases. The average velocity, $U_{jet}$ was also between $8 m/s$ to $16 m/s$ per each case. Thus, the time-dependent velocity magnitude at the inlet boundary was prescribed as follows:

\begin{equation}
\label{eq1}
\begin{aligned}
&V_{in} = U_{jet} + 0.75U_{jet}sin(2\pi f t).\\
\end{aligned}
\end{equation}

In addition to the specified flow parameters, we also considered a geometrical parameter, specifically, the ratio of the nozzle-to-jet distance to jet diameter ($H/d$ shown in Fig. \ref{FIG_schematic}), which can range $2\sim6$. This results in a total of three independent parameters, offering a broad spectrum of variation.

To comprehensively account for the impact of changes in all variables on heat transfer performance, one needs to conduct 125 simulations by exploring five different values for each parameter. Thus, to decrease the simulation scenarios without compromising representativeness, we needed to employ a suitable design of experiment (DOE) approach. These arrays, organized tables of experimental designs, are both efficient and balanced, enabling the simultaneous assessment of multiple factors. The method was proposed and developed by Genichi Taguchi \citep{taguchi1989quality, taguchi1995quality, taguchi2002mahalanobis} which is an iterative statistical method utilized for the optimization of the system whose behavior is under the influence of several parameters. In this work, the Taguchi method is employed to identify the most influential factors affecting the performance of the feed-forward neural network. 

In other words, Taguchi method is mainly intended to perform an efficient sampling procedure to reduce the number of effective configurations of the parameters in the data set. In this regard, as presented in Table \ref{tab:table_cases}, a reduced data set has been created by Taguchi method to avoid the possible issues \citep{singh2023fsnet, dhillon2020convolutional, jin2023novel} related to the presence of numerous features as inputs to the MLP.     

\begin{table}[ht]
    \centering
    \begin{filecontents*}{Cases_Phase1.csv}
Case,H/d,Fr,V
1,2,5,8
2,2,25,16
3,2,50,14
4,2,75,12
5,2,100,10
6,3,5,10
7,3,25,8
8,3,50,16
9,3,75,14
10,3,100,12
11,4,5,12
12,4,25,10
13,4,50,8
14,4,75,16
15,4,100,14
16,5,5,14
17,5,25,12
18,5,50,10
19,5,75,8
20,5,100,16
21,6,5,16
22,6,25,14
23,6,50,12
24,6,75,10
25,6,100,8
\end{filecontents*}
\pgfplotstabletypeset[
        col sep=comma,
        columns={Case,H/d,Fr,V}, 
        every head row/.style={before row=\hline,after row=\hline},
        every last row/.style={after row=\hline},
    ]{Cases_Phase1.csv} 
    \caption{Cases generated for harmonic jet by Taguchi method}
    \label{tab:table_cases}
\end{table}

As outlined in Table \ref{tab:table_cases}, each input consist of three values: $H/d$, $Fr$, and $V$. The corresponding output is the Nusselt number, which varies continuously through a cycle. To adapt this continuous data for the MLP, we face the challenge of sampling it across the cycle. Due to the irregular nature of the Nusselt variations, simplifying it into a single set of data based on frequency, amplitude, and phase, similar to the input, is not feasible.

One approach is to sample the data continuously through the cycle, generating, for example, 100 discrete data points. However, maintaining a balanced dimensionality between the input and output (label) of the MLP is crucial for optimal performance. Significant differences may prevent the network from learning complex patterns and hinder convergence.

To address this issue, we opted for the augmented MLP-FFT method. In this approach, the output Nusselt number is transformed into the Fourier space, providing sets of frequencies, amplitudes, and phase values. The data is organized based on amplitude, retaining only the top 10 amplitudes while discarding the others. This allows us to characterize the output Nusselt distribution using 10 sets of frequencies, amplitudes, and phases. With this transformation, we can effectively apply the MLP between the inputs and outputs. It is important to highlight that, for MLP validation, 21 out of the 25 cases were chosen for training, while the remaining four were employed for testing.

In Fig. \ref{FIG_loss_MLP} the training and test loss of the MLP has been presented. As seen below, the test loss decreases with the number of epochs approaching that of the train suggesting the successful performance of the MLP in capturing the underlying patterns and trends in the data. 

\begin{figure}
	\centering
	\includegraphics[scale=0.4]{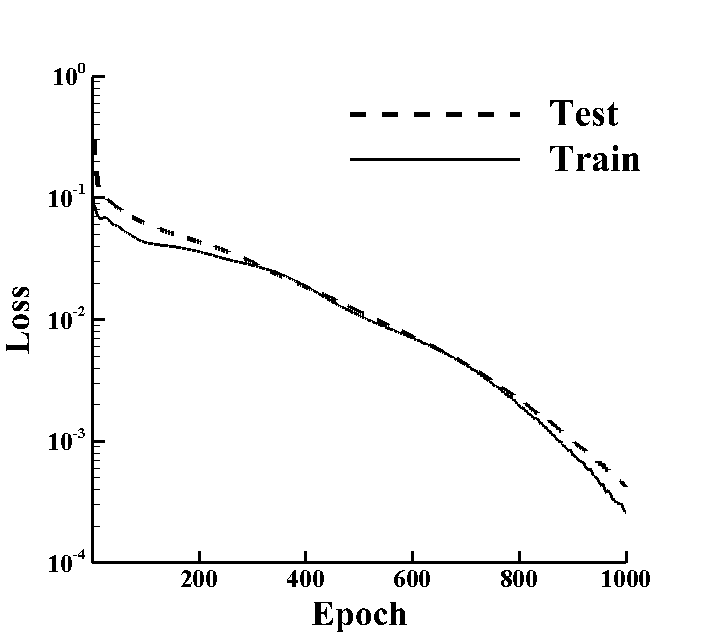}
	\caption{Train and Test loss for the MLP utilized in the constant-frequency scenario.}
	\label{FIG_loss_MLP}
\end{figure}

\begin{figure*}
    \centering
    \begin{subfigure}[b]{0.45\textwidth}
        \centering
        \includegraphics[width=0.8\textwidth]{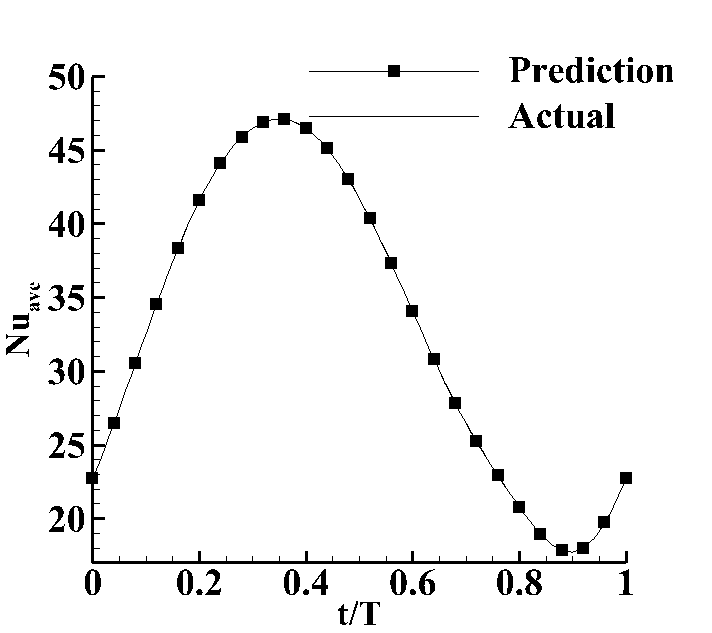}
        \caption{case21}
        \label{fig:subfig-a}
    \end{subfigure}
    \hfill
    \begin{subfigure}[b]{0.45\textwidth}
        \centering
        \includegraphics[width=0.8\textwidth]{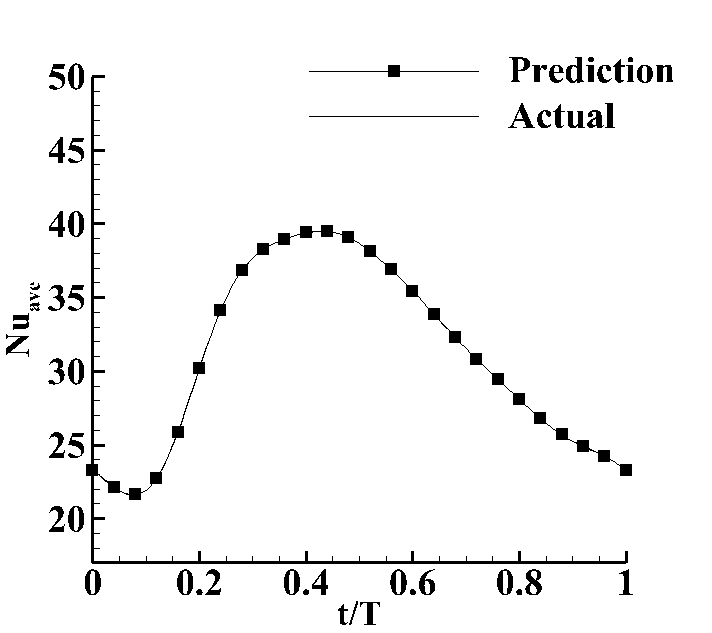}
        \caption{case22}
        \label{fig:subfig-b}
    \end{subfigure}
    \\
    \begin{subfigure}[b]{0.45\textwidth}
        \centering
        \includegraphics[width=0.8\textwidth]{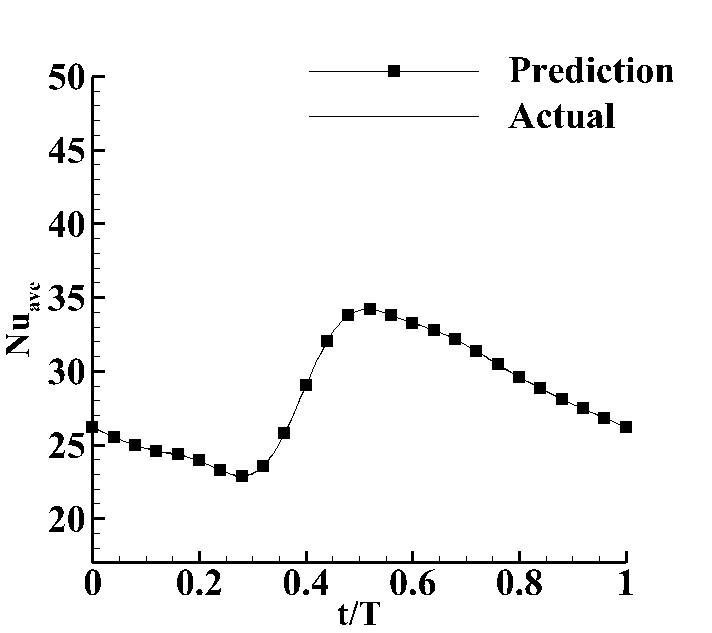}
        \caption{case23}
        \label{fig:subfig-c}
    \end{subfigure}
    \hfill
    \begin{subfigure}[b]{0.45\textwidth}
        \centering
        \includegraphics[width=0.8\textwidth]{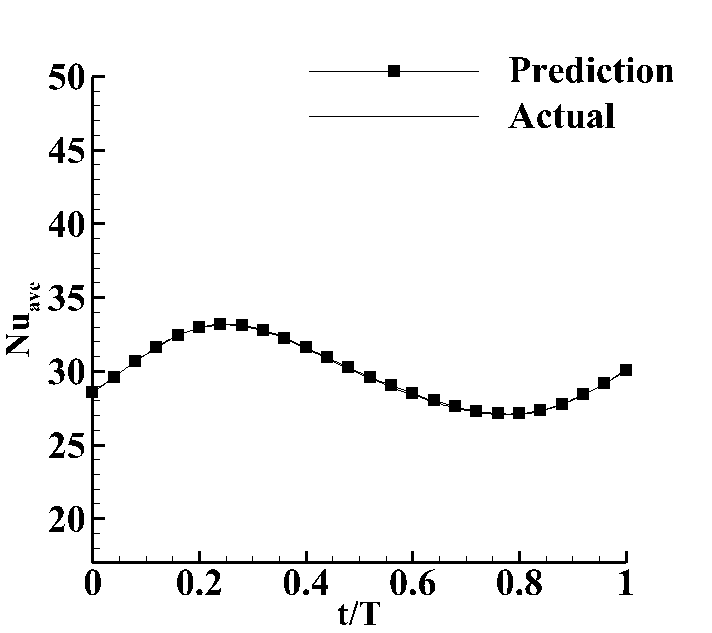}
        \caption{case24}
        \label{fig:subfig-d}
    \end{subfigure}
    \hfill
    \begin{subfigure}[b]{0.45\textwidth}
        \centering
        \includegraphics[width=0.8\textwidth]{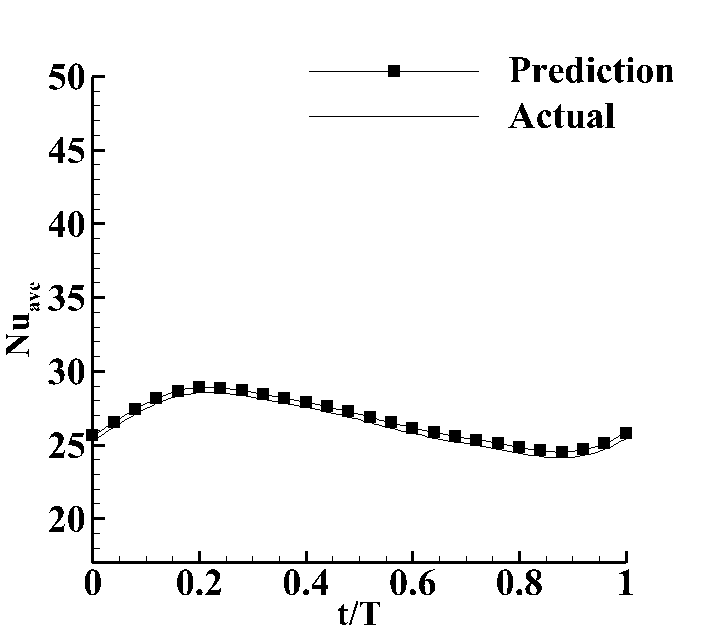}
        \caption{case25}
        \label{fig:subfig-e}
    \end{subfigure}
    \caption{Comparison between the average Nu obtained from FOM and the MLP prediction the test cases}
    \label{fig:combined}
\end{figure*}

The performance of the MLP with respect to the actual data has been presented in Fig. \ref{fig:combined}. From the figure, it is evident that the MLP has achieved a remarkably high level of precision in capturing the Nu number fluctuations. The close trends and values between the predicted results and the actual data demonstrate the effectiveness of the MLP in accurately modelling and reproducing the observed variations. Therefore, it is found that MLP can be regarded as an accurate tool to predict Nu number under single frequency inlet conditions. In other words, the strong agreement between the predicted values and the actual data supports the notion that the MLP can be considered as a robust and precise tool for Nu number prediction in such scenarios. This highlights the potential of MLP-based models as valuable tools in the field of fluid dynamics and heat transfer analysis. 

\subsection{Random-Frequency Impingement Jet}
\subsubsection{Average Nu Number Prediction}
To take the effect of varying frequency of the impingement jet into consideration, the inlet velocity was assumed to be composed of some harmonic functions with different frequencies as follows:
\begin{equation}
\label{Rndom_eq}
\begin{aligned}
&V_{in} = \sum[
U_{i} + 0.75U_{i}sin(2\pi f_{i} t))] ,\\
\end{aligned}
\end{equation}
where these values were considered for the average velocity and frequency:
\begin{equation}
\label{eq_reange}
\begin{aligned}
&U_i = \{8, 10, 12, 14, 16\} (m/s),\\
&f_i = \{5, 25, 50, 75, 100\} (Hz).\\
\end{aligned}
\end{equation}

\begin{figure}
    \centering
    \begin{subfigure}[a]{0.4\textwidth}
        \centering
        \includegraphics[width=0.8\textwidth]{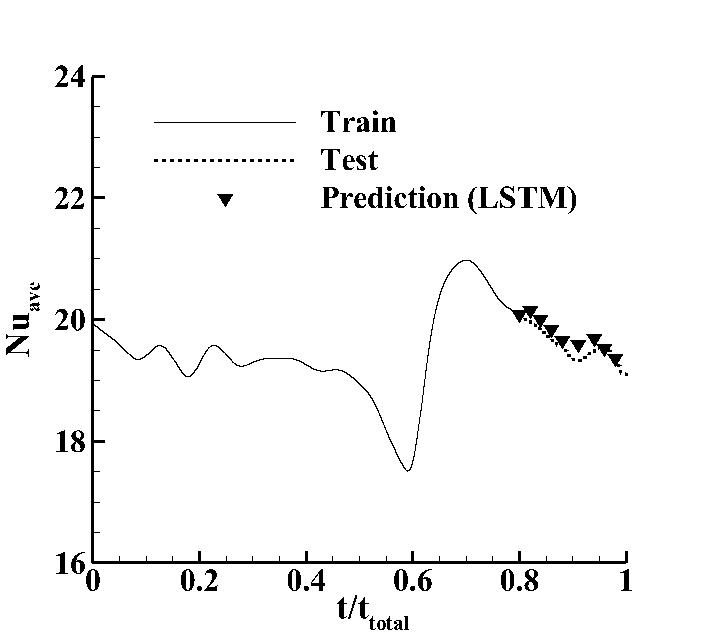}
        \caption{LSTM}
        \label{fig:subfig-a_LSTM}
    \end{subfigure}
    \hfill
    \begin{subfigure}[b]{0.4\textwidth}
        \centering
        \includegraphics[width=0.8\textwidth]{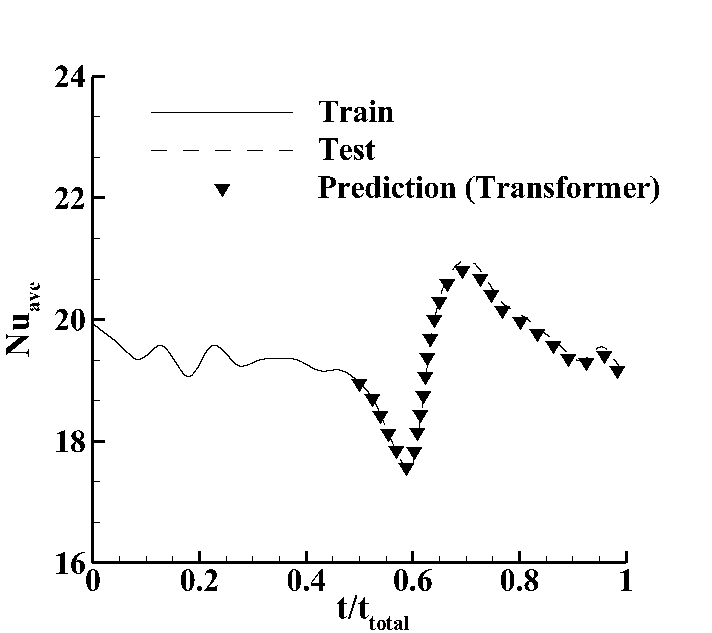}
        \caption{Transformer}
        \label{fig:subfig-b_Transformer}
    \end{subfigure}
    \\
    \caption{Comparison between the prediction made by Transformer and LSTM}
    \label{fig:LSTM_Transformer}
\end{figure}

To predict the behavior of the average Nu on the plate under these conditions, two predictive networks, namely LSTM and Transformer, were employed. The objective was to draw a comparison between their respective performances in predicting the temporal variation of the average Nu on the plate. The train and test loss for both methods has been presented in Fig. \ref{fig:LSTM_Transformer}. It should be noted that, in the case of LSTM, only the final 20 percent of the data was employed for testing, while the Transformer network, leveraging its enhanced performance attributed to the attention mechanism, was engaged for the predictions of 50 percent of the data. As is evident in Fig. \ref{fig:LSTM_Transformer}, the test loss in both of the methods demonstrates a double descent trend. As a result, the early stopping technique was implemented to avoid poor performance \citep{nakkiran2021deep, belkin2020two, loog2020brief, fei2017adaptive}. Presented in Fig. \ref{fig:LSTM_Transformer} demonstrates the average Nu number predicted by the methods. It is evident in Fig. \ref{fig:LSTM_Transformer} that the Transformer NN has been considerably more successful in predicting the future of Nu number compared to LSTM. To illustrate, the Transformer managed to predict a longer portion of the Nu number compared to LSTM. As seen in Fig.\ref{fig:LSTM_Transformer}, almost half of the Nu number variation as a time series can be accurately predicted by the transformer, while this amount falls to less than twenty percent as seen in Fig.\ref{fig:subfig-a_LSTM}. More importantly, Transformer succeeds in predicting the amount of Nu number at each time instant with higher precision with respect to the performance of LSTM. As seen in Fig.\ref{fig:subfig-b_Transformer}, the Nu numbers predicted by Transformer fall perfectly on top of the real data obtained from the numerical simulation. However, the prediction offered by LSTM shows some discrepancies concerning the real data. Therefore, it is concluded that Transformer could be regarded as a more accurate and robust tool for predicting the Nu number on the plat subjected to a varying impingement jet compared to LSTM. This is attributed to the fact that the architecture of Transformer enjoys the attention mechanism which enables the possibility for the model to prioritize the importance of each time step in the time series. Thus, the attention mechanism leads to a more robust model owing to its capability to adapt the model to different time-dependent variability.      
\begin{figure}
    \centering
    \begin{subfigure}[a]{0.4\textwidth}
        \centering
        \includegraphics[width=0.8\textwidth]{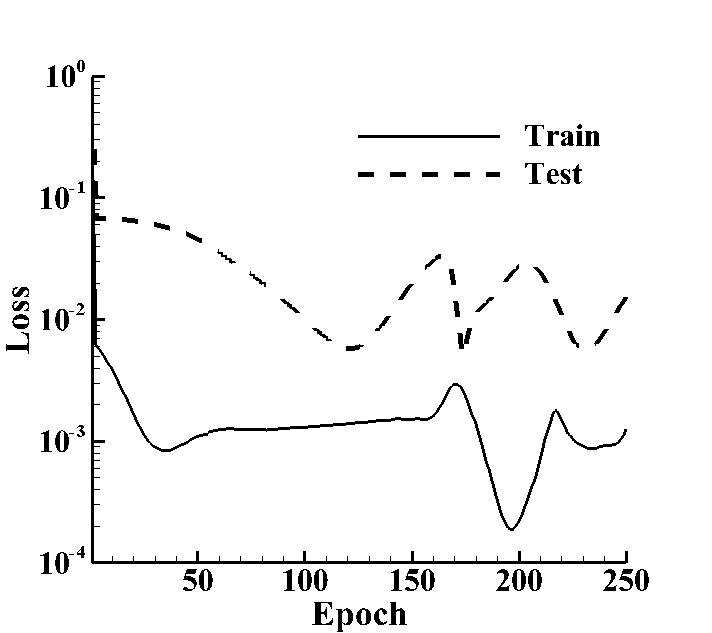}
        \caption{Loss for LSTM}
        \label{fig:subfig-a_LSTM_loss}
    \end{subfigure}
    \hfill
    \begin{subfigure}[b]{0.4\textwidth}
        \centering
        \includegraphics[width=0.8\textwidth]{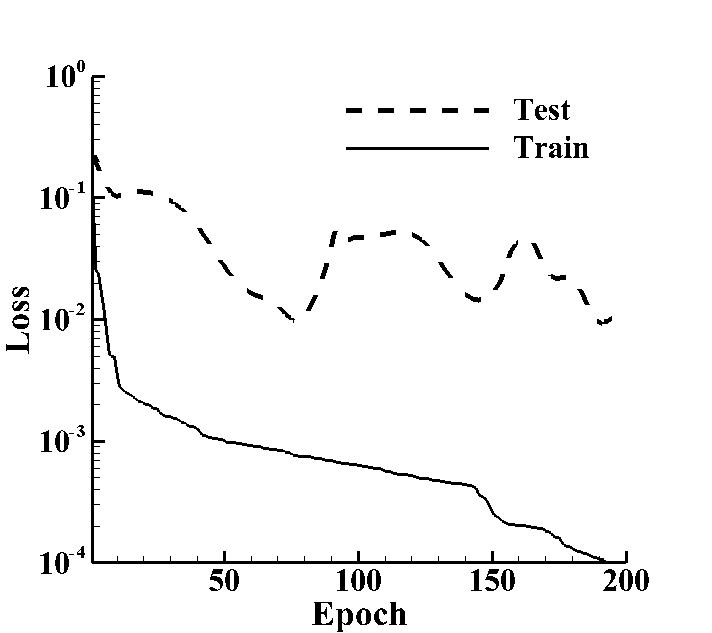}
        \caption{Loss for Transformer}
        \label{fig:subfig-b_Transformer_loss}
    \end{subfigure}
    \\
    \caption{Comparison between the train and test loss for Transformer and LSTM}
    \label{fig:LSTM_Transformer_loss}
\end{figure}

\subsubsection{Local Nu Number Prediction}
\begin{figure}
	\centering
	\includegraphics[scale=0.35]{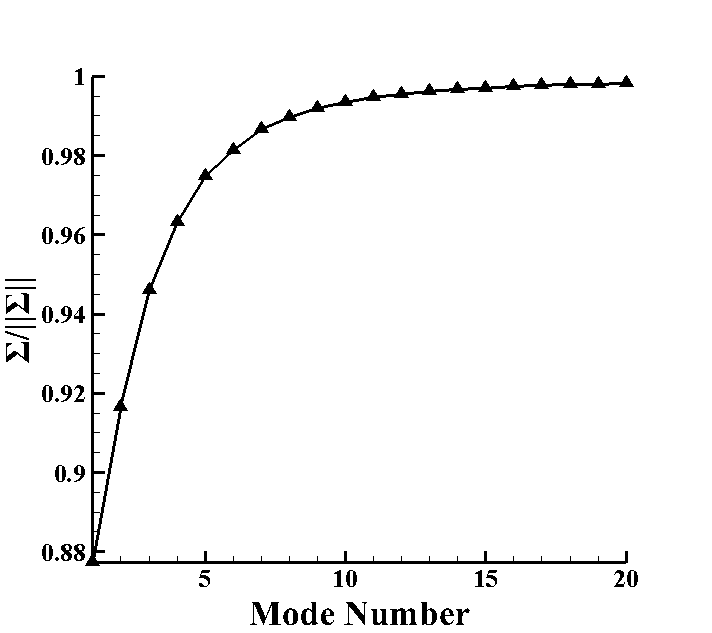}
	\caption{Cumulative Modal Energy}
	\label{cum_enrg}
\end{figure}

\begin{figure*}
	\centering
	\includegraphics[scale=0.9]{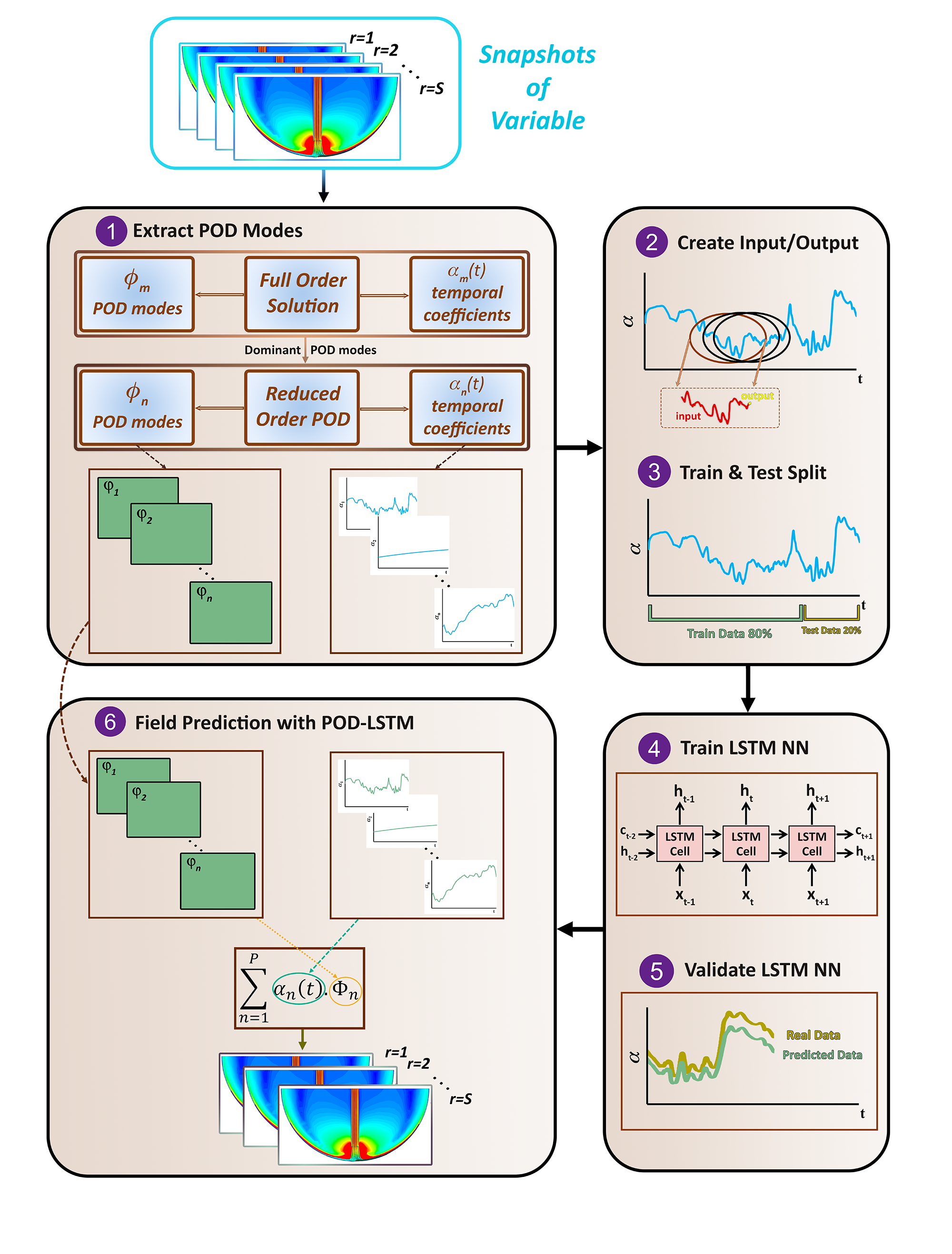}
	\caption{Structure of the POD-LSTM approach.}
	\label{FIG_POD_LSTM}
\end{figure*}

The techniques discussed so far deal with forecasting the average Nu on the plate. Despite the cardinal importance of this number in engineering applications, gaining information on local Nu is crucial in many cases. In this regard, a novel method is presented in this work to work out this problem. Fundamentally, it should be noted that treating a contour as a time series is associated with some difficulties in reaching the desired accuracy especially, due to the numerous structures being present in the data structures. To address this issue, the POD+LSTM method was utilized which can be summarized in these steps, according to Fig. \ref{FIG_POD_LSTM}:

\begin{enumerate}
  \item The snapshots of the variable, in this case, the Nusselt number on the symmetry boundary, have been gathered in a matrix. The POD analysis is then performed on this matrix to define the spatial orthogonal modes $\phi_n$ and temporal coefficients $\alpha_n$. The number of modes $N$ is selected in a way to maintain 99\% of energy. Since $N$ is much smaller than the dimension $M$ of $\phi_n \in R^{M}$, the subspace expanded by $\phi_n$ ($n=1,2,...N$) can be regarded as the low-dimensional representation of the full order model. $\alpha_n$ are then used for training in the LSTM network. 
  It should be noted that only five modes of the snapshots matrix were considered as it corresponds to 99\% of energy, as seen in Fig. \ref{cum_enrg}.
  \item The time series of the temporal coefficients $\alpha_n$ are divided into several time windows (here we take a time window of 20).
  \item The time coefficients are then divided, with the initial 80\% of data chosen for training, while the subsequent 20\% is reserved for testing.
  \item  The LSTM network is trained in this part. The idea is to teach the LSTM neural network the relationship between the input and output pairs.
  \item The validation of the LSTM network will be performed using the predicted time coefficients with that of the test datasets.
  \item At the final step, the field will be reconstructed by using the dominant POD modes (step 1) and the predicted time coefficients obtained by the LSTM network.
\end{enumerate}


\begin{figure*}
    \centering
    \begin{subfigure}[b]{0.45\textwidth}
        \centering
        \includegraphics[width=0.7\textwidth]{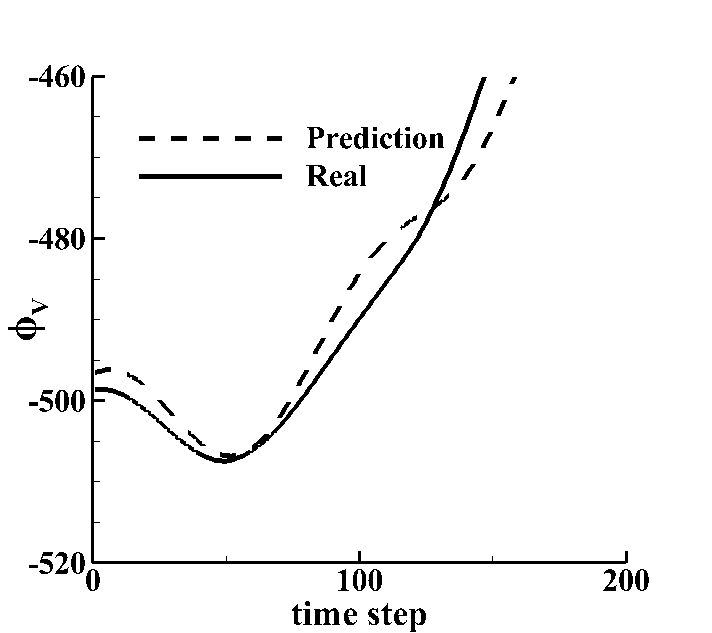}
        \caption{Mode 1}
        \label{mode1}
    \end{subfigure}
    \hfill
    \begin{subfigure}[b]{0.45\textwidth}
        \centering
        \includegraphics[width=0.7\textwidth]{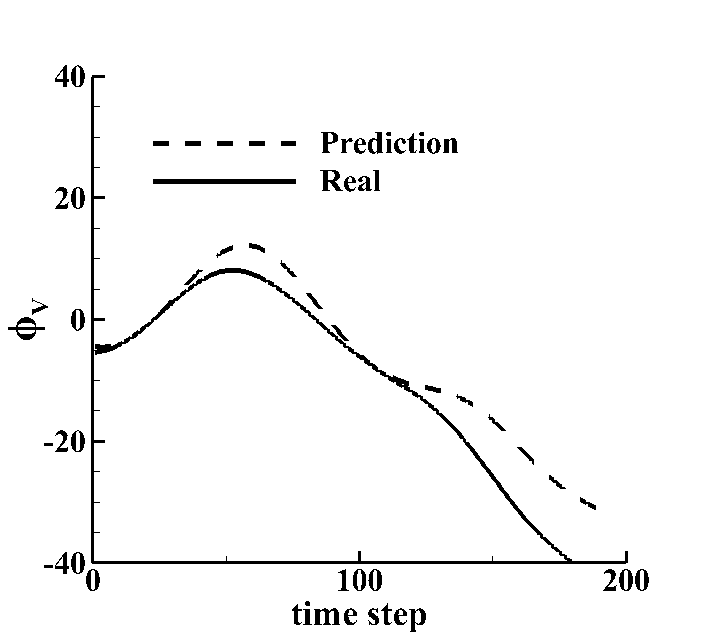}
        \caption{Mode 2}
        \label{mode2}
    \end{subfigure}
        \begin{subfigure}[b]{0.45\textwidth}
        \centering
        \includegraphics[width=0.7\textwidth]{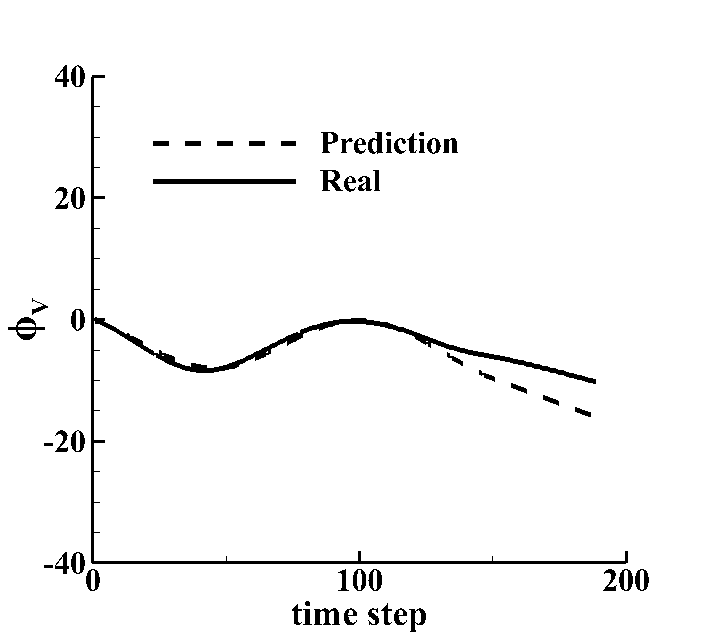}
        \caption{Mode 3}
        \label{mode3}
    \end{subfigure}
        \begin{subfigure}[b]{0.45\textwidth}
        \centering
        \includegraphics[width=0.7\textwidth]{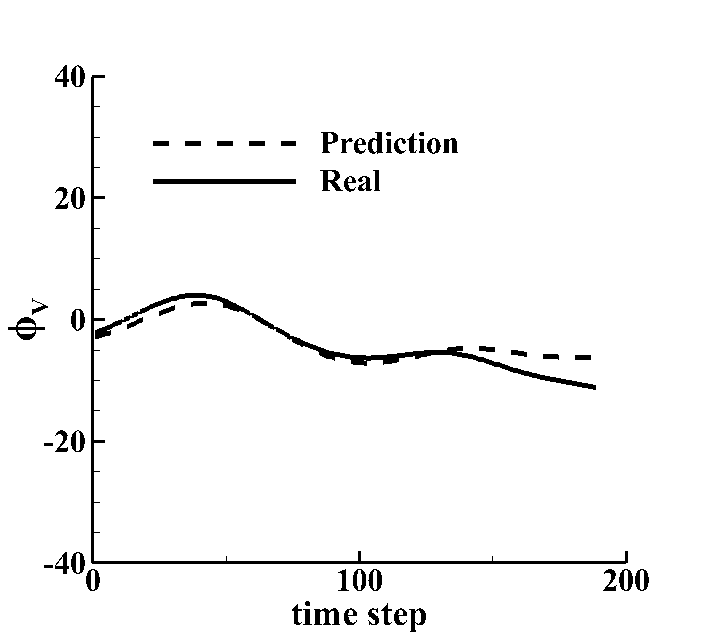}
        \caption{Mode 4}
        \label{mode4}
    \end{subfigure}
        \begin{subfigure}[b]{0.45\textwidth}
        \centering
        \includegraphics[width=0.7\textwidth]{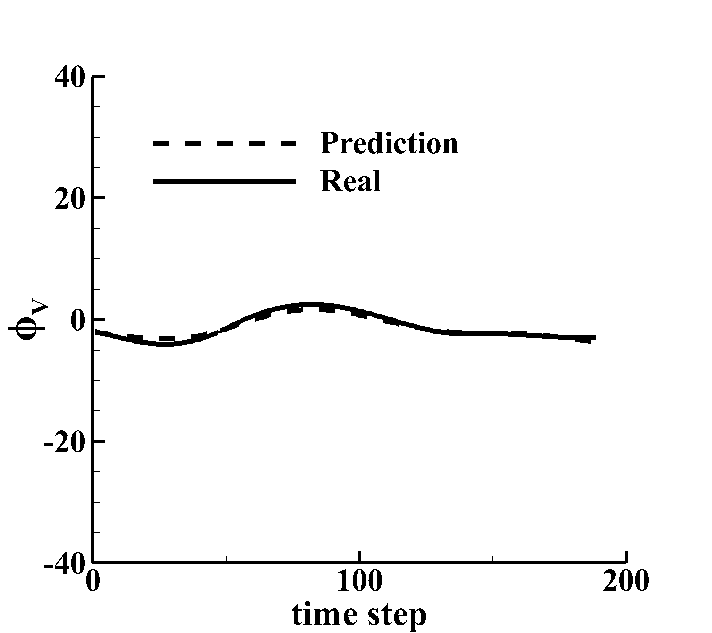}
        \caption{Mode 5}
        \label{mode5}
    \end{subfigure}
    \\
    \caption{Comparison between actual and predicted values for different time coefficients}
    \label{fig:FOM_ROM}
\end{figure*}

In the following, the viability of this method has been evaluated. The prediction of the temporal modes corresponding to the first five dominant modes has been presented in Fig. \ref{fig:FOM_ROM} which proves the successful performance of LSTM in capturing both trend and value of each temporal modes in time. Furthermore, it is essential to analyze the accumulation of deviations between the actual and predicted values through time. Notably, the most substantial deviation is observed for the second mode, reaching a magnitude of 25\% at the final time.

The impact of these deviations on the final result is completely linked to the singular values (denoted as $\Sigma$ in Eq. \ref{eq_pod3}) associated with each mode. For instance, referring to the cumulative Modal Energy curve in Fig. \ref{cum_enrg}, the effect of deviation in the first mode (with a singular value of approximately 0.88) can be more than 20 times that of the second mode (singular value of 0.04). 

Unfortunately, the \textit{black box} nature of predictive neural networks, including LSTM, prevents fine-tuning to prioritize one variable over another. This limitation is particularly pronounced in scenarios where the influence of certain modes needs to be prioritized, as demonstrated by the cumulative Modal Energy curve.

\begin{figure}
    \centering
    \begin{subfigure}[a]{0.4\textwidth}
        \centering
        \includegraphics[width=0.8\textwidth]{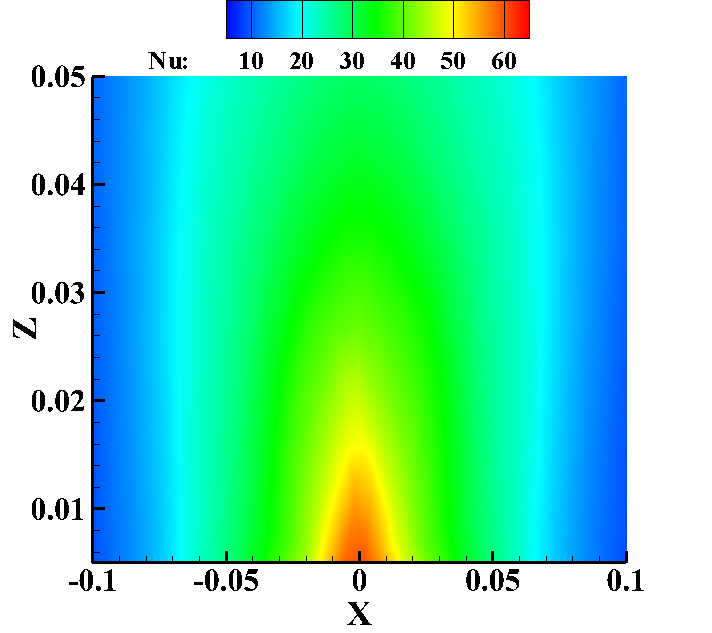}
        \caption{Full Order Solution}
        \label{FOM_pred}
    \end{subfigure}
    \hfill
    \begin{subfigure}[b]{0.4\textwidth}
        \centering
        \includegraphics[width=0.8\textwidth]{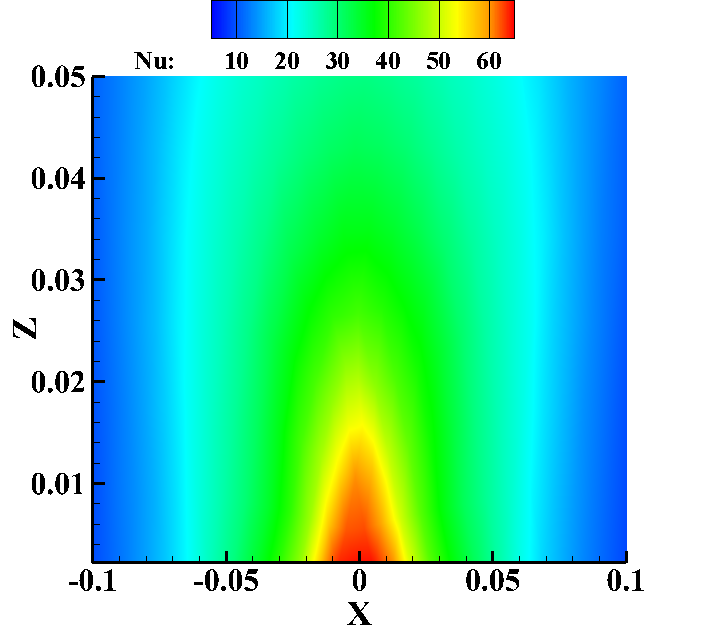}
        \caption{Predictive Surrogate Model solution}
        \label{ROM_pred}
    \end{subfigure}
    \hfill
    \begin{subfigure}[b]{0.4\textwidth}
        \centering
        \includegraphics[width=0.8\textwidth]{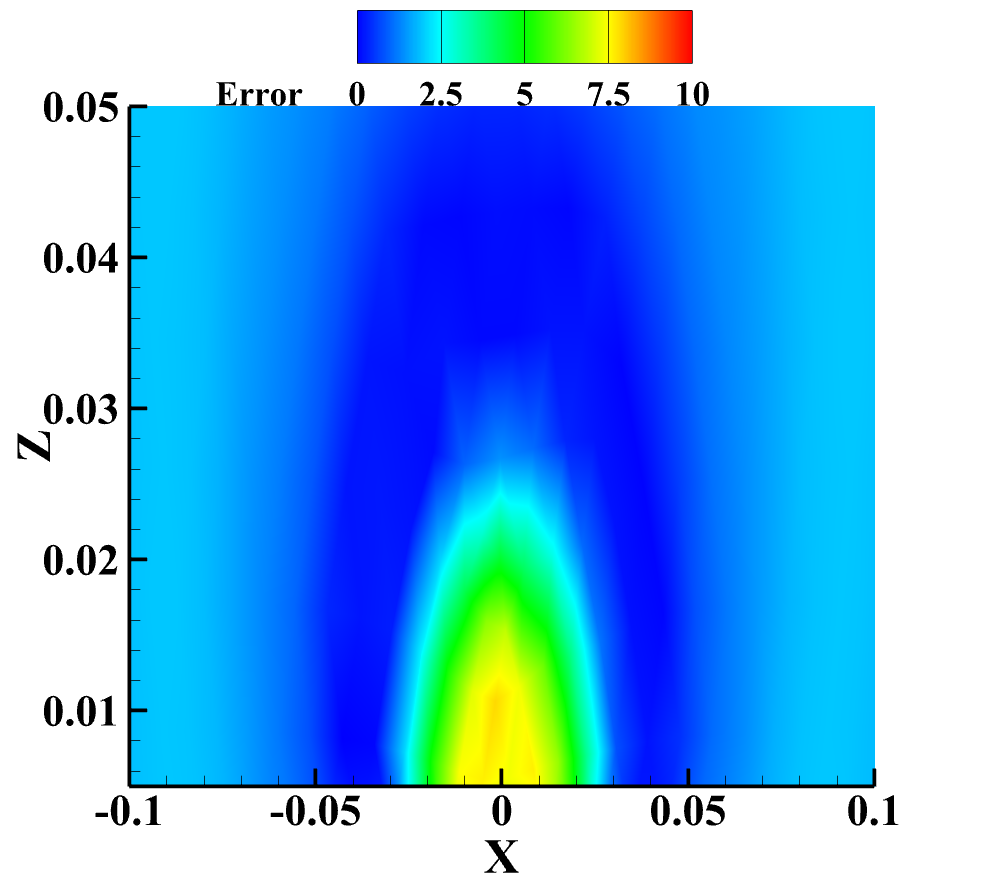}
        \caption{Error (\%)}
        \label{L2_error}
    \end{subfigure}
    \\
    \caption{Comparison between full order modelling and ROM-based prediction}
    \label{fig:FOM_ROM_contour}
\end{figure}

To visually illustrate the performance, Fig. \ref{fig:FOM_ROM_contour} compares the contour of temperature for both the FOM and the proposed predictive surrogate model (POD-LSTM) at $t=0.95t_{final}$. The method proposed in this work perfectly captures the physics of heat transfer on the plate at the consecutive time steps. This proves the practicality of POD for boosting the performance of time-series prediction methods especially, for local Nu on a plate subjected to an impingement jet. 
Additionally, the relative L2 norm error is incorporated in Fig. \ref{L2_error}. It reveals that the maximum error, approximately 8\%, is concentrated in the region directly in front of the impingement jet. And, as one moves away from the jet, the error gradually decreases, reaching as low as 0\%. This observation underscores the effectiveness of the proposed POD-LSTM model in accurately predicting temperature distribution, all over the domain.


\section{Conclusion}
In this paper, the heat transfer rate of a concave surface prone to an impinging jet is determined by various predictive surrogate models. The study systematically explores the impact of various model order reduction and deep learning techniques under different scenarios. For the constant-frequency impingement jet, the augmented MLP-FFT method demonstrates remarkable precision in capturing Nu number fluctuations. The study further extends to the random-frequency impingement jet, comparing the performance of LSTM and Transformer. The Transformer outperforms LSTM, showcasing its accuracy and robustness in predicting Nu numbers under varying impingement conditions. The Transformer demonstrates precise Nu number predictions for almost 50\% of the cycle, whereas LSTM predictions cover only 20\% of the cycle and exhibit a higher level of error, reaching a maximum of 5\%. Additionally, a novel approach combining POD and LSTM for local Nu prediction proves successful, enhancing time-series prediction for complex structures. This work underscores the potential of modern MOR and ML methods in advancing our understanding of heat transfer phenomena, which can be extended to other applications in fluid dynamics and heat transfer analysis.

\bibliographystyle{unsrtnat}

\bibliography{refs}


\end{document}